\documentclass{amsart}
\usepackage{amssymb, epic, eepic}

\newtheorem{thm}{Theorem}
\newtheorem{cor}[thm]{Corollary}
\newtheorem{lem}[thm]{Lemma}
\newtheorem{prop}[thm]{Proposition}

\newtheorem{remarkk}[thm]{Remark}
\newenvironment{remark}{\begin{remarkk} \em}{\end{remarkk}}
\newtheorem{examplee}[thm]{Example}
\newenvironment{example}{\begin{examplee} \em}{\end{examplee}}

\newcommand{\bbz}{\mathbb{Z}}
\newcommand{\bbr}{\mathbb{R}}

\title{Virtual Knot Groups}
\author{Se-Goo Kim}
\address{Department of Mathematics, Indiana University, Bloomington,
Indiana 47405}
\email{sekim@indiana.edu}


\begin{document}

\begin{abstract}
Virtual knots, defined by Kauffman, provide a natural generalization of
classical knots. Most invariants of knots extend in a natural way to give
invariants of virtual knots. In this paper we study the fundamental groups
of virtual knots and observe several new and unexpected phenomena.

In the classical setting, if the longitude of a knot is trivial in the
knot group then the group is infinite cyclic. We show that for any
classical knot group there is a virtual knot with that group and trivial
longitude. It is well known that the second homology of a classical knot
group is trivial. We provide counterexamples of this for virtual knots.

For an arbitrary group $G$, we give necessary and sufficient conditions
for the existence of a virtual knot group that maps onto $G$ with
specified behavior on the peripheral subgroup. These conditions simplify
those that arise in the classical setting.
\end{abstract}

\maketitle

\vspace*{-2em}

\section{Introduction}
\label{sec:introduction}
In 1996 Kauffman presented the theory of Virtual Knots.
In~\cite{kauffman}
he described this very natural generalization of classical knot
theory and began to develop the theory of virtual knotting.  A
significant recent accomplishment is the work of Goussarov, Polyack, and
Viro \cite{gpv} in which it is shown that
the entire theory of finite type
invariants of knots extends naturally to the
realm of virtual knots.  In this paper we explore the
fundamental groups of virtual knots, focusing on properties of
the peripheral subgroup and on the homology theory of these
groups.

The core idea of virtual knot theory is easily described.  To
each knot diagram in $\bbr^3$ one can form an associated Gauss
diagram, a simple diagram that captures all the crossing
information of the knot diagram.  Details will be given later. 
Knot diagrams determine isotopic knots if they are related by
Reidemeister moves; there are corresponding moves on Gauss
diagrams that generate an equivalence relation on the set of
Gauss diagrams.  It can be shown that if two knots determine
equivalent Gauss diagrams, they represent the same knot.

A virtual knot is defined to be an equivalence class of Gauss
diagrams.  It follows from the discussion above that every knot
determines a unique virtual knot and if two knots determine the
same virtual knot then they are in fact the same knot. 
However, not every virtual knot arises from a knot and hence
virtual knot theory offers a nontrivial extension of classical
knot theory. The work of \cite{kauffman,gpv} demonstrates that
much of
classical knot theory extends to the virtual setting.  We will
see here that while the theory of knot groups does extend in a
natural way as well, a number of new phenomena arises that
contrasts sharply with what occurs in the classical setting.

To each virtual knot $K$ there is associated a fundamental
group $\Pi_K$ along with a peripheral subgroup, generated by a
meridian and longitude, $m$ and $l$.  We observe here the
relatively simple fact that as in the classical case, $l$ and
$m$ commute.  
Classical knot groups have Wirtinger presentation of deficiency~1
and all the group
homology of dimension greater than~1 is trivial.
However, we observe that virtual knot groups
may have Wirtinger presentations of deficiency~0 and
that any Wirtinger group of deficiency~$0$ or~$1$ can be realized as
a virtual knot group. From this, we can find
examples of virtual knot groups with nontrivial second homology.

It is well known that a classical knot with trivial longitude is
the unknot. We observe that any Wirtinger group of deficiency 1
is the group of some virtual knot with trivial longitude.
As a corollary, any classical knot group is the group of a virtual knot
with trivial longitude, which gives examples of nontrivial virtual knots
with trivial longitude.

Let $G$ be a group and let $\mu$ and $\lambda$ be elements of $G$. Is
there a virtual knot $K$ and a surjective homomorphism $\rho\!:\Pi_K\to G$
such that $\rho(m)=\mu$ and $\rho(l)=\lambda$, where $m$ and $l$ are the
meridian and longitude of $K$?
In the classical knot case, Edmonds and Livingston \cite{el} answered
this for $G=S_n$, and Johnson and Livingston \cite{jl} extended this
result to general group representations.
In the virtual case, only one of their
conditions is needed, namely, that the image of longitude commutes with
a normal generator of $G$.
In the course of proving this, we examine properties of
connected sums and present an example of a nontrivial virtual knot
which is a connected sum of two trivial virtual knots.

If $g_1$ and $g_2$ are commuting elements of a group $\Pi$, then their
Pontryagin product, defined later,
is an element $\langle g_1,g_2\rangle$ in $H_2(\Pi)$.
We observe that the second homology group of a virtual knot group is
generated by the Pontryagin product of the
meridian and longitude of the virtual knot.
In the classical case, the Pontryagin product is zero since the second
homology groups of classical knots are zero.
We finally present examples of virtual knots whose groups have the second
homology $\bbz$ or $\bbz/2$, and an example of a virtual knot which has
trivial longitude and whose group is of deficiency 0.

\textsc{Acknowledgements.} The author would like to thank Chuck Livingston
for his suggestions and discussions.

\begin{figure}
\begin{center}
\input{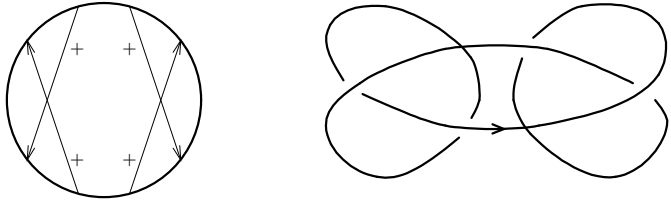}
\caption{A Gauss diagram and its corresponding virtual knot}
\label{fig:vknot}
\end{center}
\end{figure}

\section{Virtual knot groups}
\label{sec:group}
\subsection{Gauss diagrams}
Knots are usually presented by knot diagrams, that is, generic immersions
of the circle into the plane enhanced by information on overpasses and
underpasses at double points. A generic immersion of a circle into the
plane is characterized by its \emph{Gauss diagram}, which consists of the
circle together with the preimages of each double point of the immersion
connected by a \emph{chord}.
To incorporate the information on overpasses and
underpasses, the chords are oriented from the upper branch to the lower
one. Furthermore, each chord is equipped with the sign of the
corresponding double point (local writhe number). The result is called a
\emph{Gauss diagram} of the knot.
Thus Gauss diagrams can be considered as an alternative way to present
knots. However, not every Gauss diagram is indeed a Gauss diagram of
some knot. For example, see Figure~\ref{fig:vknot} and
\cite{kauffman,gpv}.

\subsection{Virtual knots}
\begin{figure}
\begin{center}
\input{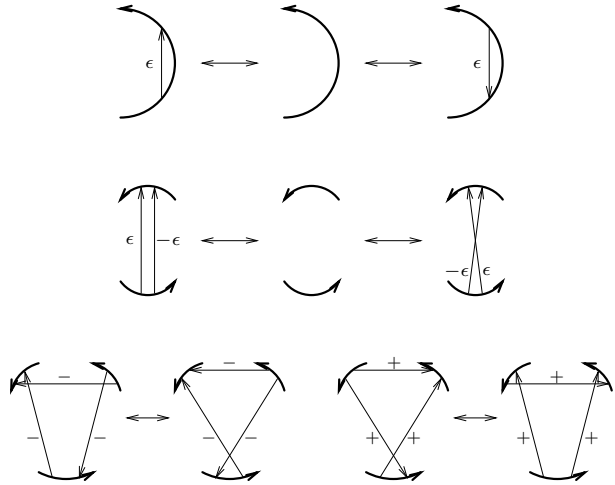}
\caption{Moves of Gauss diagrams corresponding to Reidemeister moves.}
\label{fig:reid}
\end{center}
\end{figure}

As is well-known, when a knot changes by a generic isotopy, its diagram
undergoes a sequence of Reidemeister moves. Figure~\ref{fig:reid} depicts
the counter-parts of the Reidemeister moves for Gauss diagrams. All moves
corresponding to the first and second Reidemeister moves are shown in the
top and middle rows, respectively. Though there are eight moves
corresponding to the third Reidemeister moves, it is known that only two
moves of them are necessary; they are shown in the bottom row in
Figure~\ref{fig:reid}. For details, see~\cite{gpv}. A sequence of
moves in Figure~\ref{fig:reid} is called an \emph{isotopy} for Gauss
diagrams. A \emph{virtual knot} is defined as an equivalence class of
Gauss diagrams up to isotopy. Goussarov, Polyak, and Viro~\cite{gpv}
proved that if two classical knots determine the same virtual knot,
they are isotopic in the classical sense.

\subsection{Virtual knot groups}
Kauffman~\cite{kauffman} has proved that many isotopy invariants of
knots extend naturally
to invariants of virtual knots. In particular, the notion of knot group
extends in a straightforward manner, disregarding the original
topological nature. The knot group, which is defined for classical knots
as the fundamental group of the knot complement, is extended via a
formal construction of a Wirtinger presentation (defined later).
This construction can be written down in terms
of a Gauss diagram as follows.

Let $D$ be a Gauss diagram. If we cut the circle at each arrowhead
(forgetting arrowtails), the circle of $D$ is divided into a set of
arcs. To each of these arcs there corresponds a generator of the
group. Each arrow gives rise to a relation. Suppose the sign of an
arrow is $\epsilon$, its tail lies on an arc labeled $a$, its head is
the final point of an arc labeled $b$ and the initial point of an arc
labeled $c$. Then we assign to this arrow the relation
$c=a^{-\epsilon}ba^\epsilon$.
(For simplicity, $u^{v}$ will denote $v^{-1}uv$ for any words $u$ and
$v$.) The resulting group is called the
\emph{group of the Gauss diagram}, denoted~$\Pi_D$.
One can check that it is invariant
under the Reidemeister moves of Gauss diagrams shown in
Figure~\ref{fig:reid} and that
the group of a Gauss diagram obtained from a
classical knot is the fundamental group of that knot.
Therefore, the notion of
knot group is extended to virtual knots. The \emph{group of a virtual
knot} is defined as the group of a representative Gauss diagram of the
virtual knot.
For details, see~\cite{kauffman,gpv}.

\subsection{Peripheral subgroups}
The notion of peripheral subgroup system
also extends. For the meridian, take the generator corresponding to
any of the arcs. To write down the longitude, we go along the circle
starting from this arc and write~$a^\epsilon$, when passing the head
of an arrow whose sign is $\epsilon$ and whose tail lies on the arc
labeled $a$, and finally write $t^{-p}$, where $t$ is the meridian of the
starting arc and $p$ is the number for which the result is to be in
the commutator subgroup of the virtual knot group.
The choice of such a number $p$ is possible
since the abelianization of a virtual knot group is cyclic.
One can easily check that meridian and longitude are uniquely
determined up to conjugation under the Riedemeister
moves. The \emph{peripheral subgroup} of a virtual knot is the subgroup
generated by its meridian and longitude.

\subsection{Wirtinger presentations}
Consider a group presentation of the form
$$\langle \,t_1,\ldots,t_p\mid r_1,\ldots,r_q\,\rangle,$$
where each relator $r_k,\ k=1,\ldots,q$,
is of the form $t_i^{-1}t_j^{w_k}$, for some $i$ and $j$,
$1\le i,j\le p$, and some word $w_k$ in $t_1,\ldots, t_p$.
If all $t_i$ are conjugate, such a group presentation is called a
\emph{Wirtinger presentation with respect to $t$}, where $t$ is any
element in the conjugacy class of the $t_i$.
Such an element $t$ is said to \emph{normally generate the group},
or the group is said to be \emph{of weight 1}.
A group is called \emph{Wirtinger} if it has a Wirtinger presentation.
Since the first homology group of a group is the abelianization of the
group, the first homology group of a Wirtinger group is infinite cyclic
$\bbz$.

\subsection{Deficiency of a group presentation}
For a Gauss diagram $D$, it is clear by definition that its group
$\Pi_D$ has a Wirtinger presentation with the number of generators equal
to the number of relators. The \emph{deficiency} of a group presentation is
the number of generators minus the number of relators.
The group $\Pi_D$ (hence any virtual knot group)
then has a Wirtinger presentation of nonnegative deficiency.

\begin{prop}
Any virtual knot group has a Wirtinger presentation of
deficiency 0 or 1 and hence its second homology group is cyclic.
\end{prop}

\begin{proof}
Let $\Pi$ have a Wirtinger presentation of nonnegative
deficiency $d$. Then the group
$\Pi$ has a presentation with $n+d$ generators $t_i$ and $n$ relators
$r_j$. Define a
2-dimensional CW complex $X^2_\Pi$ as follows: the 1-skeleton of
$X^2_\Pi$ is a one-point union of $n+d$ circles, each of which
represents each generator~$t_i$, and $n$ 2-disks are
attached to the 1-skeleton along the relators $r_j$. By the Van Kampen
theorem, $\Pi$ is the fundamental group of $X^2_\Pi$.
Thus, $H_1(X^2_\Pi)=\bbz$ and rank$(H_1(X^2_\Pi))=1$.
The cellular chain complex of $X^2_\Pi$ is
$\ldots\to 0\to \bbz^n \stackrel{\partial_2}{\to} \bbz^{n+d}
\stackrel{\partial_1}{\to}\bbz$, where $\partial_1$ is
a zero map, and thus
$d\le$ rank$(\mathrm{coker}\,\partial_2)=$ rank$(H_1(X^2_\Pi))=1$.
In conclusion, the group $\Pi_D$ of a Gauss diagram $D$ (hence any
virtual knot group) has a Wirtinger presentation of deficiency 0 or 1.

By rank arguments similar to the one above, $H_2(X^2_\Pi)$
is zero and infinite cyclic if the deficiency is 1 and 0, respectively.
Continuing attaching higher cells to $X^2_\Pi$
we get an Eilenberg-MacLane space $K(\Pi,1)$ whose 2-skeleton is $X^2_\Pi$.
Thus, the group homology $H_2(\Pi)$, which is defined as
$H_2(K(\Pi,1))$, is cyclic since it is a quotient of $H_2(X^2_{\Pi})$.
\end{proof}

\section{Realization}
\label{sec:realization}
In this section, we show that any Wirtinger group $\Pi$ of
deficiency 0 or 1 can be realized as a virtual knot group.
In the classical case, if $H_2(\Pi)=0$, Kervaire \cite{kervaire} showed
that $\Pi\cong\pi_1(S^n-S^{n-2})$ for
some smooth $(n-2)$-sphere $S^{n-2}$ in $S^n$, $n\ge 5$.
However, Kervaire's theorem does not apply to the case $n=3$
and there are Wirtinger groups of
deficiency~1 which are not classical knot groups.
Also, if $H_2(\Pi)\ne 0$, then $\Pi$ can never be realized as a
classical knot group. For details, see \cite{fox2,kervaire}.

\subsection{Cyclic and realizable Wirtinger presentations}
A \emph{cyclic} Wirtinger presentation is a
Wirtinger presentation of the form
$\langle\, t_1,\ldots, t_n\mid r_1,\ldots, r_n\,\rangle$, where the $j$-th
relator $r_j$ is of the form $t_{j+1}^{-1}t_j^{w_j}$ for
$j=1,\ldots,n$ (mod $n$), and some word~$w_j$ in $t_1,\ldots, t_n$.
A cyclic Wirtinger presentation can be transformed into
a special kind of cyclic Wirtinger presentation, called a
\emph{realizable} Wirtinger presentation,
with the property that each $w_j$ is a one-letter word $t_k^\epsilon$
for some $k=1,\ldots,n$ and $\epsilon=\pm 1$,
by introducing more generators and relators.
For example, if $w_1=t_3t_7^{-1}t_5$ and
$t_2=t_1^{w_1}$, then we introduce two more
generators, say, $t'_1,t'_2$, remove the relator $r_1$ and add three
more relators $t_1^{\prime -1}t_1^{t_3}$,
$t_2^{\prime -1}(t'_1)^{t^{-1}_7}$, and
$t_2^{-1}(t'_2)^{t_5}$. It is easy to see
that the new group presentation presents the same group.
Note that both cyclic and realizable Wirtinger presentations have
deficiency~$0$.

\begin{lem}
\label{lem:transform}
A Wirtinger presentation of deficiency 0 or 1 can be transformed to a
realizable Wirtinger presentation.
\end{lem}
\begin{proof}
Let $\Pi$ be a Wirtinger presentation of deficiency 0 or 1, say,
$\Pi=\langle\, t_1,\ldots,t_n\mid r_1,\ldots,r_m\,\rangle$, where $m=n$ or
$n-1$. If $m=n-1$, by doubling the relator $r_m$,
we may assume $m=n$. By the arguments prior to this lemma,
it suffices to construct a cyclic Wirtinger presentation from $\Pi$.

Let $P_\Pi$ be a graph with $n$ vertices 
$\{v_1,\ldots,v_n\}$ and $n$ edges corresponding to relators in the way
that an edge has end vertices $v_i$ and $v_j$ if and only if there is a
relator of the form $t_i^{-1}t_j^w$. Such an edge is denoted by
$e^i_j$. Since all $t_i$ are conjugate, the graph $P_\Pi$ is
connected. If two edges $e^i_j$ and $e^j_k$ meet at a vertex~$v_j$,
then we have two relators $t_i^{-1}t_j^{w_1}$ and
$t_j^{-1}t_k^{w_2}$, or $t_i=t_j^{w_1}$ and $t_j=t_k^{w_2}$. This implies
$t_i=w_1^{-1}w_2^{-1}t_kw_2w_1=(w_2w_1)^{-1}t_kw_2w_1$. We now remove
the relator $t_i^{-1}t_j^{w_1}$ and add
$t_i^{-1}t_k^{w_2w_1}$ to get a new presentation. It is
obvious that the new presentation presents the same group. This
operation corresponds to an operation on the graph $P_\Pi$ of deleting
an edge $e^i_j$ and adding an edge $e^i_k$.

Recall that a cycle of a graph is a simply closed path on the graph.
Using the above operation on the graph $P_\Pi$
we will construct a cycle from $P_\Pi$. Since
$P_\Pi$ has betti number $1$, it has one and only one cycle~$C$.
We will use induction on the length $l$ of $C$. If $l=n$, then
$C=P_\Pi$ and thus $P_\Pi$ is a cycle.
Suppose $l<n$. Then there is a vertex $v_i$ which is not on $C$. Since
$P_\Pi$ is connected, there is a path from $v_i$ to a vertex of
$C$. On this path, there is an edge $e^j_k$ such that $v_j$ is
not in $C$ but $v_k$ is in $C$. By the operation described above, we
can construct $P'_\Pi$ with a cycle $C'$ containing all vertices of
$C$ and $v_j$. Since the length of $C'$ is greater than the length of $C$,
we keep doing this process to get a cycle.

The corresponding operations on the Wirtinger presentation of
$\Pi$ give a cyclic Wirtinger presentation. This completes the proof.
\end{proof}

\begin{thm}
\label{thm:realization}
Any Wirtinger presentation of deficiency $0$ or $1$ can be realized as
a virtual knot group.
\end{thm}

\begin{proof}
By Lemma~\ref{lem:transform}, we can assume that a realizable Wirtinger
presentation with $n$ generators~$t_i$ and $n$ relators~$r_i$ is given.
We start with a circle. Choose $n$
points on the circle, which divide the circle into $n$ arcs
labeled $t_1,\ldots,t_n$, successively counterclockwise.
For each $i=1,\ldots,n$, if
$r_i=t_{i+1}^{-1}t_i^{t_j^\epsilon}$, attach an oriented chord with
sign $\epsilon$ from a point on the arc labeled
$t_j$ to the point dividing $t_i$ and
$t_{i+1}$. By definition, the group of this Gauss diagram
has the given presentation.
\end{proof}

\begin{figure}
\begin{center}
\input{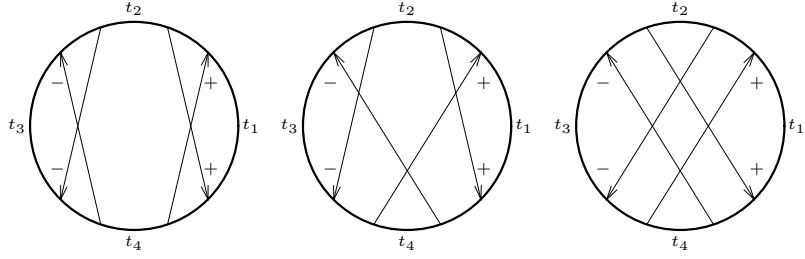}
\caption{Gauss diagrams constructed from a group presentation
$\langle\, t_1,t_2,t_3,t_4\mid
t_2^{-1}t_1^{t_4},\ t_3^{-1}t_2^{t_4^{-1}},\
t_4^{-1}t_3^{t_2^{-1}},\ t_1^{-1}t_4^{t_2}\,\rangle$}
\label{fig:real}
\end{center}
\end{figure}

\begin{example} 
The Gauss diagrams in Figure~\ref{fig:real} are
constructed from a Wirtinger presentation
$\langle\, t_1,t_2,t_3,t_4\mid
t_2^{-1}t_1^{t_4},\ t_3^{-1}t_2^{t_4^{-1}},\
t_4^{-1}t_3^{t_2^{-1}},\ t_1^{-1}t_4^{t_2}\,\rangle$.
This presentation can be transformed to
$\langle\, t_2, t_4\mid t_2t_4t_2=t_4t_2t_4\,\rangle$ that is a
presentation of the trefoil group. Note also that the longitude is
trivial. The first and second diagrams have bracket polynomial
$-A^{-6}-A^{-4}+A^{-2}+3+A^2-A^4-A^6$ while the third has unit
bracket polynomial~$1$.
Thus, there are virtual knots with the same peripheral
structure but different bracket polynomials. The third is also an example
of nontrivial virtual knot with unit bracket polynomial.
For the notion of the bracket polynomial, see~\cite{kauffman}.
\end{example}

\subsection{On peripheral subgroups}
Let $D$ be a Gauss diagram with group $\Pi_D=\langle\, t_1,\ldots,t_n \mid
r_1,\ldots,r_n\,\rangle$, where $r_i=t^{-1}_{i+1}t_i^{w_i}$ for each
$i=1,\ldots,n$ (mod $n$) and some word $w_i$ in $t_1,\ldots,t_n$. Then the
longitude of $D$ is $l=w_1\cdots w_n t_1^{-p}$, where
$p$ is chosen so that $l$ is in the commutator subgroup of $\Pi_D$.

\begin{prop}
\label{prop:peripheral}
The longitude of $D$ commutes with the meridian $t_1$ and hence the
peripheral subgroup is abelian.
\end{prop}
\begin{proof}
Iterating all relators $r_1,\ldots, r_n$, we have
$t_1=w_n^{-1}\cdots w_1^{-1}t_1w_1\cdots w_n$, or $t_1$ commutes with the
longitude. This then shows that the peripheral subgroup of $\Pi_D$ is an
abelian group generated by the meridian and longitude.
\end{proof}

\begin{prop}
\label{prop:trivial}
If a Wirtinger presentation has deficiency $1$ and if $\lambda$ is a
commutator element commuting with a normal generator $t$,
then it is the group of a virtual knot with longitude $\lambda$.
In particular, a Wirtinger presentation of deficiency $1$
can be realized as the group of a virtual knot with trivial longitude.
\end{prop}

\begin{proof}
Suppose that a group $\Pi$ has a Wirtinger presentation of deficiency 1.
Then adding a redundant relator $r_n=t_1^{-1}w_1\cdots w_{n-1} t_n
w_{n-1}^{-1}\cdots w_1^{-1}$ to the presentation gives a Wirtinger
presentation of deficiency~0 presenting the same group.
We can easily see that the longitude obtained from the new presentation
is trivial.
Moreover, in this case, if $\lambda$ is a commutator element
commuting with $t_1$
in $\Pi$, adding $r_n=t_1^{-1}\lambda^{-1} w_1\cdots w_{n-1} t_n
w_{n-1}^{-1}\cdots w_1^{-1}\lambda$ to the original presentation, we have
longitude $\lambda$ for the virtual knot obtained from the presentation.
\end{proof}

\begin{cor}
Any classical knot group is the group of a virtual knot with
trivial longitude.
\end{cor}

\begin{cor}
There is a nontrivial virtual knot with trivial longitude.
\end{cor}

We will see a partial converse of Proposition~\ref{prop:trivial} in
Section~\ref{sec:pontryagin}. An example of a Wirtinger presentation of
deficiency $0$ having trivial longitude will be presented in
Section~\ref{sec:examples}.

\section{Peripherally specified homomorphs}
\label{sec:peripherally}
Let $G$ be a group and let $\mu$ and $\lambda$ be elements of $G$. Is
there a virtual knot $K$ and a surjective homomorphism $\rho\!:\Pi_K\to G$
such that $\rho(m)=\mu$ and $\rho(l)=\lambda$, where $m$ and $l$ are the
meridian and longitude of~$K$? In the classical knot case
Edmonds and Livingston~\cite{el} answered this for $G=S_n$,
and Johnson and Livingston \cite{jl} extended this result to general group
representations.
In this section, we observe that only one of their conditions
is needed in the virtual setting,
namely, that the image of longitude commutes with
a normal generator of $G$.

\subsection{Realizable elements}
An initial observation is that $G$ must be finitely generated and of
weight one. Fix $G$ and a normal generator $\mu$ with these properties.
If, for given $\mu$ and $\lambda$ in $G$, a virtual knot $K$ and
representation $\rho$ as above exist, we say $\lambda$ is
\emph{realizable}. Let $\Lambda$ denote the set of realizable elements.
Johnson and Livingston have proved the following in the classical knot case.

\begin{thm}
\emph{\cite{jl}} \label{thm:jl1}
The set of realizable by classical knot groups is a nonempty subgroup of $G$.
\end{thm}
\begin{thm}
\emph{\cite{jl}} \label{thm:jl2}
$\lambda$ is realizable by a classical knot
if and only if $\lambda\in G''\cap Z(\mu)$,
$\langle\mu,\lambda\rangle=0\in H_2(G)$, and
$\{\mu,\lambda\}=0\in H_3(G/G')/p_\ast(H_3(G))$.
\end{thm}
In the above statement $G''$ denotes the second commutator subgroup of
$G$ and $Z(\mu)$ is the centralizer of~$\mu$. $\langle\ ,\ \rangle$
denotes the Pontryagin product (which will be defined later),
which maps pairs of commuting elements
of~$G$ to~$H_2(G)$. $\{\ ,\ \}$ is a map which sends pairs of elements
of~$G$ satisfying the first and second conditions to
$H_3(G/G')/p_\ast(H_3(G))$. $p\!:G\to G/G'$ is the natural projection.
$G'$ denotes the commutator subgroup of~$G$.

The counter-parts of the above two theorems for virtual knots are the
following.

\begin{thm}
\label{thm:realizable}
$\Lambda$ is a nonempty subgroup of $G$.
\end{thm}

\begin{thm}
\label{thm:peripheral}
$\Lambda=G'\cap Z(\mu)$, \emph{i.e.}, a commutator element~$\lambda$
in~$G$ is realizable if and only if
it commutes with the normal generator $\mu$.
\end{thm}
In order to prove these, we introduce a based connected sum of two Gauss
diagrams.

\subsection{Connected sum}
For $i=1,2$, let $D_i$ be a Gauss diagram and let $p_i$ be a point
of the circle of $D_i$ meeting no chords. We assume that each circle of
$D_i$ has a counter-clockwise orientation. The \emph{connected sum of
$D_1$ and $D_2$ based at $(p_1,p_2)$} is defined in usual manner: cut small
intervals around $p_i$ not intersecting any chords and attach the ends of
the intervals according to the orientations.

\subsection{Group presentation of a connected sum}
Suppose that the groups of two Gauss diagrams $D_1$ and $D_2$ are
$\Pi_{D_1}=\langle\,t_1,\ldots,t_n\mid
t_2^{-1}t_1^{u_1},\ldots,t_1^{-1}t_n^{u_n} \,\rangle$
and $\Pi_{D_2}=\langle\,s_1,\ldots,s_m\mid
s_2^{-1}s_1^{v_1},\ldots,s_1^{-1}s_m^{v_m} \,\rangle$, respectively,
where $u_i$ are words in $t_j$ and $v_i$ are words in $s_j$.
Suppose that $p_1$ is on the arc represented by $t_1$ and $p_2$ is on the
arc represented by $s_1$. Then the group of the
connected sum $D$ of $D_1$ and
$D_2$ at $(p_1,p_2)$ has a presentation
$$\langle\, t_1,\ldots,t_n,s_1,\ldots,s_m\mid
t_2^{-1}t_1^{u'_1},\ldots,s_1^{-1}t_n^{u'_n},\
s_2^{-1}s_1^{v'_1},\ldots,t_1^{-1}s_m^{v'_m}
\,\rangle,$$
where $u_i'$ is obtained from $u_i$ in such a way that each $t_1$ in $u_i$
is replaced with $s_1$ if it is read off from a chord with arrow tail
lying on the arc between the end point of the arc $t_n$
and the point $p_1$, and $v_i'$ is defined in the similar way.

\begin{figure}
\begin{center}
\input{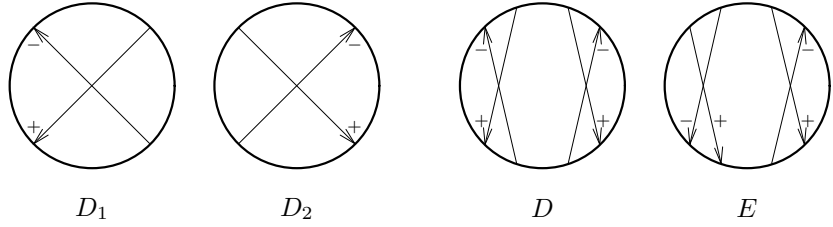}
\caption{Connected sums of two trivial virtual knots $D_1$ and $D_2$}
\label{fig:sum}
\end{center}
\end{figure}

\begin{example}
A striking example is a nontrivial virtual knot that is a connected sum
of two trivial virtual knots. The Gauss diagram $D$ in
Figure~\ref{fig:sum} has trefoil group while $E$ can be
transformed into a trivial Gauss diagram by Reidemeister moves.
This also says that an unbased connected sum of two virtual knots
is not uniquely defined.
\end{example}

\subsection{Representation of a connected sum}
Let $\rho_i\!:\Pi_{D_i}\to G\ (i=1,2)$ be representations
with $\rho_1(t_1)=\rho_2(s_1)=\mu$ and
$\rho_i(l_i)=\lambda_i$, where $l_i$ is the longitude to $D_i$.
Define a representation $\rho$ of the group of the connected sum $D$
of $D_1$ and $D_2$ at $(p_1,p_2)$ to $G$ by
$\rho(t_i)=\rho_1(t_i),\ i=1,\ldots,n$ and
$\rho(s_j)=\rho_2(s_j),\ j=1,\ldots,m$.
Since $\rho_1(t_1)=\rho_2(s_1)$, the representation $\rho$ is well-defined
and it satisfies that $\rho(m)=\mu$ and
$\rho(l)=\lambda_1\lambda_2$, where $m=t_1$ and $l$ are the meridian and
longitude of~$D$.
This is the key to proving Theorem~\ref{thm:realizable}.

\begin{proof}[Proof of Theorem~\ref{thm:realizable}.]
If $\lambda_1$ and $\lambda_2$ are realizable, then a connected sum can be
used to show $\lambda_1\lambda_2$ is realizable. To realize
$\lambda_1^{-1}$, use the Gauss diagram used to realize~$\lambda_1$
with orientation and all signs of chords reversed. That the set is nonempty
follows from Theorem~\ref{thm:jl1} since a knot group is a virtual knot
group. It follows from a method of modifying group representations
similar to the ones given in Section~\ref{sec:realization} as well.
\end{proof}

\begin{remark}
Consider a knot $K$ in $S^3$ with the meridian $m$ and longitude $l$ and a
representation $\rho\!:\pi_1(S^3\setminus K)\to G$ such that $\rho(m)=\mu$
and $\rho(l)=\lambda$ for some~$\lambda$ in $G$. The connected sum
of $K$ and its mirror image $\bar{K}$ gives a representation
$\rho\!:\pi_1(S^3\setminus K\#\bar{K})\to G$ mapping the longitude to the
identity element in~$G$. In particular, the identity element is realizable
by a classical knot group. Note that this follows from
Theorem~\ref{thm:jl1} as well.
\end{remark}

\begin{proof}[Proof of Theorem~\ref{thm:peripheral}.]
It has been proved in Section~\ref{sec:realization} that
$\Lambda\subset G'\cap Z(\mu)$. To show the converse,
let $\lambda$ be an element in $G'\cap Z(\mu)$.
Let $K$ be a knot in $S^3$ and let $\rho\!:\pi_1(S^3\setminus K)\to G$ be
a representation with $\rho(l)=1$ as seen in the remark prior to this
proof. The knot group $\Pi=\pi_1(S^3\setminus K)$ has a Wirtinger
presentation $\langle\, t_1,\ldots,t_n\mid t_2^{-1}t_1^{w_1},\ldots,
t_n^{-1}t_{n-1}^{w_{n-1}}, t_1^{-1}t_n^{w_n}\,\rangle$,
where the last relator is redundant,
the longitude of $K$ is $l=t_1^{-p}w_1\cdots w_n$, $\rho(t_1)=\mu$ and
$\rho(l)=1$.

Since $\rho$ is onto and hence the restriction of $\rho$ to $\Pi'$
is onto $G'$, there is an element~$u$ in $\Pi'$ with
$\rho(u)=\lambda$. Let $\Xi=\langle\, t_1,\ldots,t_n\mid
t_2^{-1}t_1^{w_1},\ldots,t_n^{-1}t_{n-1}^{w_{n-1}},
t_1^{-1}t_n^{w_n u}\,\rangle$. Then $\Xi$ is a quotient group
of $\Pi$ with an extra relator $t_1^{-1}t_n^{w_n u}$.
Thus $\Xi$ has a Wirtinger presentation of deficiency 0.
Since $lu\in\Xi'$, the group $\Xi$
can be realized as the group of a virtual knot with longitude $lu$.
Define $\bar{\rho}\!:\Xi\to G$ by $\bar{\rho}(t_i)=\rho(t_i)$,
$i=1,\ldots,n$. Since $\bar{\rho}(t_1^{-1}t_n^{w_n u})=
(\bar{\rho}t_1)^{-1}((\bar{\rho}t_n)^{\bar{\rho}w_n})^{\bar{\rho}u}=
\mu^{-1}\mu^\lambda=1$, $\bar\rho$ is a well-defined surjective
homomorphism and the longitude maps to $\bar\rho(l)\bar\rho(u)=\lambda$.
This implies that $\lambda$ is realizable.
\end{proof}

\section{The Pontryagin product and the second homology group of
a virtual knot}
\label{sec:pontryagin}
\subsection{Pontryagin product}
If $g_1$ and $g_2$ are commuting elements of a group $\Pi$, then their
\emph{Pontryagin product}
is an element $\langle g_1,g_2\rangle \in H_2(\Pi)$
defined as follows:
Because $[g_1,g_2]=1$ there is a homomorphism $\phi\!:\bbz\times\bbz\to\Pi$
given by $\phi(1,0)=g_1$ and $\phi(0,1)=g_2$.
Now $H_2(\bbz\times\bbz)=\bbz$
with generator $z$ described as the image of $1\otimes 1$ under the cross
product isomorphism $H_1(\bbz)\otimes H_1(\bbz)\to H_2(\bbz\times\bbz)$
arising in the K\"{u}nneth formula. Set $\langle
g_1,g_2\rangle=\phi_\ast(z)$. It can be interpreted geometrically as follows:
Since $[g_1,g_2]=1$ in $\Pi$, $[g_1,g_2]$ is nullhomotopic in the
Eilenberg-MacLane space $K(\Pi,1)$. Thus there is a map $\varphi$ from a
torus $T^2$ to $K(\Pi,1)$ sending a meridian and a longitude to paths
representing $g_1$ and $g_2$. The Pontryagin product $\langle
g_1,g_2\rangle$ is $\varphi_\ast[T^2]$ in $H_2(K(\Pi,1))$,
where $[T^2]$ is the generator of
$H_2(T^2)$. A reference is \cite{brown}.

\subsection{A CW-complex obtained from a realizable Wirtinger presentation}
Let $D$ be a Gauss diagram with group
$\Pi=\langle\, t_1,\ldots,t_n\mid t_2^{-1}t_1 ^{t_{i_1}^{\epsilon_1}},
\ldots,
t_1^{-1}t_n^{t_{i_n}^{\epsilon_n}}\,\rangle$.
We define a 2-dimensional CW complex whose fundamental group is $\Pi$.
Consider a torus $T^2$ as a square identified linearly the top side with
the bottom side and the left side with the right side.
We divide the square into $n$ rectangles $R_1,\ldots,R_n$
by $n-1$ vertical lines and label the vertical lines $t_1,\ldots, t_n$
successively from the left to the right with orientation from the bottom to
the top. A CW complex $X_D$ is defined as a quotient of $T^2$ as
follows:
For $j=1,\ldots, n$, the top side of each rectangle $R_j$
is oriented from the left to
the right and identified linearly with the left vertical line
$t_{i_j}$ of the rectangle $R_{i_j}$ in an
oriented manner according to the sign of $\epsilon_j$.
In fact, if $D$ represented a classical knot $K$ in $S^3$, $X_D$ is the
2-skeleton of a CW complex of $S^3\setminus K$.
It can be easily checked that $\Pi$ is the fundamental group of $X_D$.

\begin{figure}
\begin{center}
\input{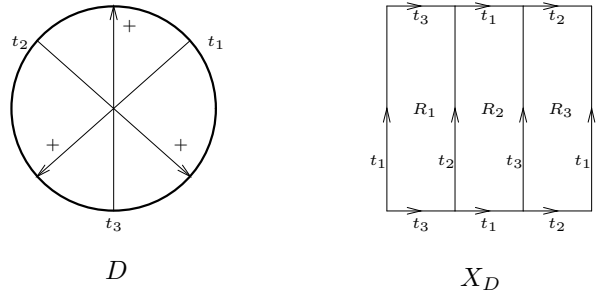}
\caption{Gauss diagram $D$ and space $X_D$}
\label{fig:space}
\end{center}
\end{figure}

\subsection{A description on the second homology of a group}
As briefly described in Section~\ref{sec:group},
there is a description on the second homology of a group $\Pi$. If
$Y$ is a connected CW-complex with $\pi_1(Y)\cong\Pi$ and $\Sigma_2(Y)$
denotes the subgroup of $H_2(Y)$ generated by all singular 2-cycles
represented by maps of a 2-sphere into $Y$, then
$H_2(\Pi)=H_2(Y)/\Sigma_2(Y)$. Note $\Sigma_2(Y)=\rho(\pi_2(Y))$,
where $\rho\!:\pi_2(Y)\to H_2(Y)$ is the Hurewicz homomorphism.
Compare \cite{kervaire,bms}.

\begin{thm}
The Pontryagin product of the meridian and longitude of a virtual knot
generates the second homology of its group.
\end{thm}

\begin{proof}
Let $D$ be a representative Gauss diagram of a virtual knot.
We first compute the second homology of $X_D$. Consider the cellular
chain complex of $X_D$
$$0 \to C_2 \stackrel{\partial_2}{\to} C_1\stackrel{\partial_1}{\to}
C_0\to 0,$$
where $C_2$, $C_1$, and $C_0$ are the free abelian groups generated by the
rectangles $R$'s, the circles $t$'s, and the vertex, respectively.
The map $\partial_1$ is a zero map and the map $\partial_2$ is represented
by a matrix
{\scriptsize$\left(\begin{array}{rrrr}
1 & 0 & \cdots & -1 \\
-1 & 1 & \cdots & 0 \\
0 & -1 & \cdots & 0 \\
\vdots & \vdots & & \vdots \\
0 & 0 & \cdots & 1
\end{array}\right)$}
where $\ker\partial_2$ is generated by the element
$[R_1+\cdots+R_n]$. Let $q$ be the quotient map from $T^2$ to $X_D$. Then
$q_\ast[T^2]$ is the element $[R_1+\cdots+R_n]$ in $H_2(X_D)$,
where $[T^2]$ denotes the generator of $H_2(T^2)$.

Observe that a meridian of the Gauss diagram $D$ is $t_1$ and its
longitude can be exactly read off from the top side of the square of $X_D$.
Since the meridian and the longitude of $T^2$ map to $t_1$ and the top
side of the
square via the quotient map $q\!:T^2\to X_D$, respectively,
$q_\ast[T^2]$ in $H_2(X_D)$ is the Pontryagin product
of the meridian and longitude of the diagram $D$ and hence it generates
$H_2(X_D)$ and $H_2(\Pi)=H_2(X_D)/\Sigma_2(X_D)$.
\end{proof}

\begin{cor}
\label{cor:pontryagin}
If a virtual knot has a nonzero second homology group,
the Pontryagin product of its meridian and longitude is not zero.
\end{cor}

The following is a partial converse of Proposition~\ref{prop:trivial}
because the second homology of a Wirtinger presentation of deficiency 1
is zero.

\begin{cor}
If a group $\Pi$ is the group of a virtual knot with trivial longitude,
then $H_2(\Pi)=0$.
\end{cor}

We conclude this section posing a question:
Can any Wirtinger presentation of deficiency 0 with trivial
second homology be realized as the group of a virtual knot with trivial
longitude?

\section{Examples}
\label{sec:examples}
\subsection{Examples of virtual knot groups with nontrivial second homology}
Edmonds and Livingston \cite{el} showed that
the Pontryagin product of the meridian and longitude of a
classical knot group represented to a symmetric group~$S_n$ is zero,
and Johnson and Livingston \cite{jl} extended this result to a general
group representation. However, there are virtual knots whose second
homology is nonzero and hence whose Pontryagin
products of meridians and longitudes are not zero. By
Theorem~\ref{thm:realization} and Corollary~\ref{cor:pontryagin}
it suffices to find examples of Wirtinger presentations of
deficiency~$0$ with nontrivial second homology.

\subsection{Infinite cyclic second homology}
Gordon \cite{gordon} gave a family of
Wirtinger presentations of deficiency~$0$ whose
second homology groups are infinite cyclic $\bbz$.
The groups are defined by, for $k\ge 2$,
$$\langle\, t,x,y\mid t^{-1}x^kt=x^{k+1},\ t^{-1}y^kt=y^{k+1},\
t^{-1}xyt=xy \,\rangle$$
Setting $z=tx$ and $w=ty$ gives
the Wirtinger presentation
$$\langle\, t,z,w\mid z=t^{(t^{-1}z)^{-k}},\
w=t^{(t^{-1}w)^{-k}},\ t=t^{w^{-1}tz^{-1}} \,\rangle$$
of deficiency 0.

\subsection{Order two second homology}
Brunner, Mayland, and Simon \cite{bms} gave a
Wirtinger group of deficiency $0$ whose second homology group is cyclic of
order 2. The group is defined by
$$\langle\, a,b\mid b=a^{-1}b^2ab^{-2}a,\
b=[ba^{-1},a^{-1}b]^{-1}b[ba^{-1},a^{-1}b]\,\rangle$$
which is a Wirtinger presentation of deficiency 0.

I do not know any examples of Wirtinger group presentations
of deficiency 0 whose second homology groups are cyclic
of finite order other than 2.

\subsection{An example of a virtual knot with trivial longitude whose
group has deficiency~$0$}
The group $\Pi$ with presentation $\langle\, x,a \mid xa^2=ax,\ a^2x=xa\,
\rangle$ considered by Fox in~\cite{fox2}, Example 12, has vanishing
second homology group by Kervaire's Theorem in~\cite{kervaire}
because it is the group of some 2-knot.
However, this group does not have deficiency 1 (it has therefore
deficiency 0) because its Alexander ideal, which is easily computed to be
$(3,1+t)$, fails to be principal.
(For the notion of the Alexander ideal, see \cite{fox1}.)
By letting $y=ax$, the presentation is transformed into a Wirtinger
presentation $\langle\,x,y\mid y=x^{y^{-1}x},\ x=y^{xy^{-1}}\,\rangle$.
Its longitude is $y^{-1}x^2y^{-1}$.
To see that the longitude is the identity in $\Pi$,
consider the infinite cyclic cover $\tilde{X}_\Pi$ of the space~$X_\Pi$.
The cover $\tilde{X}_\Pi$ has a group presentation $\langle\,
x_i,i\in\bbz\mid x_i^{-1}x_{i+1}^2,\ x_j^{-2}x_{j+1}\,\rangle$, which is
isomorphic to $\bbz/3\bbz$. The longitude $y^{-1}x^2y^{-1}$ lifts to
$x_0x_1=x_0^3=1$ which is trivial in $\bbz/3\bbz$ and hence trival in
$\Pi$ because $\pi_1(\tilde{X}_\Pi)$ is the commutator subgroup of $\Pi$.
Therefore, there is a virtual knot with trivial longitude whose group has
deficiency 0 and a trivial second homology.


\begin{thebibliography}{99}
\bibitem{brown} K. Brown, Cohomology of Groups, GTM 87, Springer-Verlag,
New York, 1982.
\bibitem{bms} A. M. Brunner, E. J. Mayland, Jr., and J. Simon, \emph{Knot
groups in $S^4$ with nontrivial homology}, Pacific J. Math. \textbf{103}
(1982), no. 2, 315--324.
\bibitem{el} A. Edmonds and C. Livingston, \emph{Symmetric representations
of knot groups}, Topology Appl. \textbf{18} (1984), 281--314.
\bibitem{fox1} R. H. Fox, \emph{Free differential calculus II}, Ann. of
Math. \textbf{59} (1954), 196--210.
\bibitem{fox2} R. H. Fox, \emph{A quick trip through knot theory}, Topology
of 3-manifolds, M. K. Fort, Jr. ed., Prentice-Hall, 1962, 120--167.
\bibitem{gordon} C. Mc. A. Gordon, \emph{Homology of groups of surfaces
in the 4-sphere}, Math. Proc. Cambridge Philos. Soc. \textbf{89} (1981),
no. 1, 113--117.
\bibitem{gpv} M. Goussarov, M. Polyak and O. Viro, \emph{Finite type
    invariants of classical and virtual knots}, preprint (October 1998
  -- math.GT/9810073).
\bibitem{jl} D. J. Johnson and C. Livingston, \emph{Peripherally
specified homomorphs of knot groups}, Trans. Amer. Math. Soc. \textbf{311}
(1989), no. 1, 135--146.
\bibitem{kauffman} L. Kauffman, \emph{Virtual Knot Theory}, preprint
(November 1998 -- math.GT/9811028).
\bibitem{kervaire} M. A. Kervaire, \emph{On higher dimensional knots},
Differential and Combinatorial Topology, A Symposium in Honor of Marston
Morse, S. S. Cairns ed., Princeton University Press, Princeton, New Jersey,
1965, 105--119.
\end{thebibliography}
\end{document}